\documentstyle[12pt,amssymb,amscd,psfig]{amsart}
\newtheorem{thm}{Theorem}
\newtheorem{lemma}[thm]{Lemma}
\newtheorem{definition}[thm]{Definition}
\newtheorem{proposition}[thm]{Proposition}

\theoremstyle{definition}

\theoremstyle{remark}

\newcommand{\co}{\colon\thinspace}    
\newcommand{\np}{\newpage}            

\title[Dwyer's filtration and $4$-manifolds]{Dwyer's filtration and topology 
of ${\mathbf 4}$-manifolds}
\author{Vyacheslav S. Krushkal}
\thanks{Partially supported by NSF grant DMS 00-72722}
\address{Department of Mathematics, University of Virginia, Charlottesville, VA 22904}
\email{krushkal\char 64 virginia.edu}

\begin{document}

\begin{abstract}  
Topological $4$-dimensional surgery is conjectured to fail, in general,
for free fundamental groups. M. Freedman and P. Teichner have shown 
that surgery problems with an arbitrary fundamental group have a solution, 
provided they satisfy a certain condition on Dwyer's filtration on 
second homology. We give a new geometric proof of this result, and analyze
its relation to the canonical surgery problems.
\end{abstract}

\maketitle

The lower central series of the fundamental group of a space $X$
is closely related to the Dwyer's \cite{D} filtration ${\phi}_k(X)$
of the second homology $H_2(X; {\mathbb{Z}})$.
It is well known \cite{FQ} that, for any $k>1$, the canonical $4$-dimensional 
surgery problems may be arranged to have the kernel represented by
a submanifold $M$ which is ${\pi}_1$-null, and satisfies 
$H_2(M)={\phi}_k(M)$. It is conjectured \cite{F} that these canonical
problems do not have a solution. On the other hand, Freedman and 
Teichner showed in \cite{FT} that if the surgery kernel is ${\pi}_1$-null
and its second homology lies in the $\omega$-term of the filtration,
$H_2(M)={\phi}_{\omega}(M)$, then the surgery problem has a solution.

This theorem was a development of the earlier result (\cite{FQ}, Chapter 6)
that a surgery problem can be solved if the kernel $M$ is ${\pi}_1$-null
and $H_2(M)$ is {\em spherical}. In the present paper we give a new,
geometric, proof of the theorem of Freedman -- Teichner. 
We show that a surgery problem with the kernel 
$M$ in their setup, i.e. ${\pi}_1$-null and with 
$H_2(M)={\phi}_{\omega}(M)$, in a $4$-manifold $N$, can be 
reduced to the ${\pi}_1$-null {\em spherical} case, in the 
same manifold $N$. The proof is based on the idea of splitting
of capped gropes. This technique has been useful in solving a 
number of other problems in $4$-manifold topology, see \cite{K}, \cite{KQ}.

It is interesting to note that the proof goes through if 
$H_2(M)={\phi}_k(M)$ for $k\geq k_0$, where $k_0$ is a constant depending 
on the inclusion $M\subset N$. This integer $k_0$ can be easily 
read off from the data as the number of group elements in ${\pi}_1 N$,
represented by the double point loops of the Whitney disks
for $M$ in $N$. As is common in the subject, the aforementioned 
canonical surgery problems may be chosen to satisfy 
$H_2(M)={\phi}_{k_0-1}(M)$, just missing the requirement for 
the theorem above.

Section \ref{gropes} gives a brief overview of the splitting operation 
on gropes. The background material on the lower central series and Dwyer's
filtration is presented in section \ref{dwyer}. The main theorem 
is stated and proved in section \ref{surgery}.

\section{Gropes and splitting} \label{gropes}

This is a brief summary of terminology and notations; for a more detailed
exposition the reader is referred to \cite{FQ}, \cite{FT}, \cite{K}, 
\cite{KQ}. Observe that the notion of a {\em grope of class $k$}
which is needed in this paper is different from the {\em symmetric gropes}
discussed in \cite{FQ}, \cite{KQ}. In particular, a symmetric grope
of height $n$ has class $2^n$. To be precise, we recall the definition:

\vspace{.2cm}

\begin{definition} \sl
A {\bf grope} is a special pair ($2$-complex, circle). A grope has
a class $k=2,3,\ldots$. A grope of class $2$ is a compact oriented
surface $\Sigma$ with a single boundary component. 
A $k$-grope is defined inductively as follows: let $\{ {\alpha}_i, 
{\beta}_i\}$ be a standard symplectic basis of circles for $\Sigma$.
For any positive integers $p_i, q_i$ with $p_i+q_i=k$, a $k$-grope is 
formed by gluing a $p_i$-grope to each ${\alpha}_i$ and a $q_i$-grope to 
each ${\beta}_i$. 
\end{definition}

\vspace{.2cm}

The {\em tips} of a grope $g$ is a symplectic basis of circles in its top 
stage surfaces. They freely generate ${\pi}_1 g$.
A model {\em capped grope} $g^c$ is obtained from a grope $g$ by attaching
disks to its tips.
The grope $g$ is then called the {\em body} of $g^c$.
Finally, a capped grope in a $4$-manifold $M$ is an immersion
$g^c\longrightarrow M$, where only intersections among the caps are
allowed (so the body is embedded, and is disjoint from the interiors
of the caps.) Each intersection point between the caps carries
an element of ${\pi}_1 M$. It is determined by the double point loop
of the intersection.

Here we use the terminology of \cite{K}. In particular,
$g$ denotes a grope (the underlying $2$-complex), while the capital letter
$G$ indicates the use of its untwisted $4$-dimensional thickening. 
The operations that will be used in the proof are {\em contraction}, 
sometimes also referred to as 
{\em symmetric surgery}, and {\em pushoff}, which are described in 
detail in \cite[section 2.3]{FQ}. The following lemma (suitably
formulated grope splitting) is a central
ingredient in the proof of the main theorem in section \ref{surgery}.
For more applications of grope splitting see \cite{K}, \cite{KQ}.

\vspace{.2cm}

\begin{lemma}[Grope splitting] \label{splitting} \sl
Let $(g^c,{\gamma})$ be a capped grope in $M^4$. 
Then, given a regular neighborhood 
$N$ of $g^c$ in $M$, there is a capped grope $(g^c_{\rm split},{\gamma})
\subset N$, such that each cap of $g^c_{\rm split}$ has double points which 
represent at most one group element in ${\pi}_1 M$, and each body surface, above
the first stage, of $g^c_{\rm split}$ has genus $1$.
\end{lemma}

\vspace{.2cm}

{\bf Proof.} First assume that $N$ is the untwisted thickening of $g^c$, $N=G^c$,
and moreover let $g^c$ be a model capped grope (without double points).
Let $C$, $D$ be a dual pair of its caps, and let ${\alpha}$ be an arc in 
$C$ with endpoints on the boundary of $C$. 
(In our applications, $\alpha$ will be chosen to separate intersection 
points of $C$ corresponding to different group elements.) Recall that 
the untwisted thickening $N$ of $g^c$ is defined as the thickening in ${\mathbb{R}}^3$, 
times the interval $I$. We consider the $3$-dimensional thickening,
and surger the top-stage surface of $g$, which is capped by $C$ and $D$,
along the arc $\alpha$. The cap $C$ is divided by ${\alpha}$ into two 
disks $C'$, $C''$ which serve as the caps for the new grope; their dual 
caps $D'$, $D''$ are formed by parallel copies of $D$. This operation 
increases the genus of this top-stage surface by $1$.
We described this operation for a model capped grope; splitting of a 
capped grope with double points is defined as an obvious generalization.

Continue the proof of lemma \ref{splitting} by dividing each cap 
$C$ by arcs $\{ {\alpha} \}$, so that each component of $C\smallsetminus
\cup{\alpha}$ has double points representing just one group element, 
and splitting $g^c$ along all these arcs. The crucial observation is that
any future application of this technique preserves the progress achieved
up to date: the parallel copies $D'$, $D''$ of $D$ as above inherit the 
collection of the group elements carried by $D$. We apply the same operation 
to the surfaces in the $(h-1)$-st stage of the grope, separating each top 
stage surface by arcs into genus $1$ pieces. This procedure is performed
inductively, descending to the first stage of $g^c$. For example, if 
originally each body surface of a $k$-grope $g$ had genus one, and each 
cap carried $n$ group elements, then after this complete splitting procedure 
the first stage surface will have genus $n^k$. 
\qed

\vspace{.2cm}

\section{Dwyer's filtration} \label{dwyer}

\vspace{.2cm}

In this section we recall basic facts about the lower central series
and Dwyer'r filtration, and their geometric reformulation in terms
of gropes. For proofs of the propositions, see \cite{FT}. 
Recall that the lower central series of a group $H$ is defined by
${H}^1={H}$, ${H}^k=[{H}, {H}^{k-1}]$ for $k\geq 1$,
and ${H}^{\omega}=\cap_{k\in{\mathbb{N}}} {H}^k$. The following
proposition provides a geometric reformulation:

\vspace{.2cm}

\begin{proposition} \sl
A loop $\gamma$ in a space $X$ lies in ${\pi}_1(X)^k$
if and only if $\gamma$ bounds a map of some $k$-grope in $X$.
\end{proposition}

\vspace{.2cm}

Clearly, a loop $\gamma$ is in ${\pi}_1(X)^{\omega}$ iff 
for each finite $k$, $\gamma$ bounds a map of a $k$-grope in $X$.
The {\em Dwyer's subspace} ${\phi}_k(X)\subset H_2(X;{\mathbb{Z}})$ is 
defined as the kernel of the composition

$$ H_2(X)\longrightarrow H_2(K({\pi}_1 X,1))=H_2({\pi}_1 X))\longrightarrow
H_2({\pi}_1 (X)/{\pi}_1(X)^{k-1}). $$

\vspace{.5cm}

\begin{proposition} \sl
Dwyer's subspace ${\phi}_k(X)$ of $H_2(X)$ coincides with the
subset of homology classes represented by maps of closed
$k$-gropes into $X$.
\end{proposition}

\vspace{.2cm}

Here a {\em closed} $k$-{\em grope} is a $2$-complex obtained by replacing
a $2$-cell in $S^2$ with a $k$-grope. Note again that a homology class is in
the $\omega$-term of the Dwyer's filtration if and only if for each $k\geq 2$ 
it is represented by a map of a closed $k$-grope into $X$.

\np

\section{Surgery theorem for ${\phi}_{\omega}$.} \label{surgery}
Before formulating the main theorem (Theorem 1.1 in \cite{FT}), recall the 
setting for surgery. Let $N$ be a compact topological $4$-manifold, possibly 
with boundary. Suppose $f\co N\longrightarrow X$ is a degree $1$ normal map from
$N$ to a Poincar\'{e} complex $X$. Following the
higher-dimensional arguments, it is possible to find a map normally
bordant to $f$ which is a ${\pi}_1$-isomorphism, and such that
the kernel 

$$ K=ker(H_2(N; {\mathbb{Z}}[{\pi}_1 X])\longrightarrow
H_2(X; {\mathbb{Z}}[{\pi}_1 X]))$$

is a free ${\mathbb{Z}}[{\pi}_1 X]$-module. Suppose Wall's obstruction
vanishes, so there is a preferred basis for the kernel $K$ in which the
intersection form is hyperbolic. Then we say that $M\subset N$ 
{\em represents the
surgery kernel} if $M$ is ${\pi}_1$-null, $H_2(M)$ is free and
$$ H_2(M)\otimes_{\mathbb{Z}} {\mathbb{Z}}[{\pi}_1 X]\longrightarrow
H_2(N;{\mathbb{Z}}[{\pi}_1 X])$$
 
maps isomorphically onto $K$. Here we assume that $M\subset int(N)$ is a 
(compact) codimension $0$ submanifold. Recall that $M$ is ${\pi}_1$-null 
means that the inclusion induces the trivial
map ${\pi}_1 M\longrightarrow {\pi}_1 N$.

\vspace{.2cm}

\begin{thm} 
\sl Suppose a standard surgery kernel is represented by
$M\subset N$ which is ${\pi}_1$-null and satisfies 
${\phi}_{\omega}(M)=H_2(M)$. Then there is a 
normal bordism from $f\co N\longrightarrow X$ to a simple homotopy
equivalence $f'\co N'\longrightarrow X$.
\end{thm}

\vspace{.2cm}

More precisely, the proof shows that there is an integer $m$ depending on the 
inclusion $M\subset N$ so that the theorem still holds if $H_2(M)\subset{\phi}_m(M)$.

\vspace{.2cm}

{\em Proof.} Let ${\gamma}_1,\ldots, {\gamma}_k$ be loops in $M$
representing generators of ${\pi}_1 M$. Since $M$ is ${\pi}_1$-null,
there are null-homotopies ${\Delta}_1,\ldots,{\Delta}_k$ in $N$,
$\partial{\Delta}_i={\gamma}_i$. Let $f_1,\ldots, f_m\in {\pi}_1 N$ 
be the group elements represented by $M\cup{\Delta}_1\cup\ldots\cup{\Delta}_k$.
These are given by the double point loops of the null-homotopies
$\Delta$, and by the intersections ${\Delta}\cap M$. (Each component
of the intersection determines a group element by starting at the
basepoint in $M$, following a path in $M$ to ${\Delta}\cap M$ and returning
via $\Delta$, avoiding its double points. This is well-defined since 
$M$ is ${\pi}_1$-null.) 

The integer $m$ is ``universal'' for the loops in $M$: given any loop 
$\gamma$ in $M$, there is a nullhomotopy for $\gamma$ in $N$ giving rise
to at most $m$ group elements, since $M$ is ${\pi}_1$-null in
$M\cup{\Delta}$. (Represent $\gamma$ as a composition of the generators
$\{ {\gamma}_i\}$, so it bounds parallel copies of the singular
disks $\{ {\Delta}_i\}$.)

Let $G=\{ G_i\}$ be $(m+1)$-gropes representing 
standard free generators of $H_2(M)$ (corresponding to the
preferred basis of $H_2(M)\otimes_{\mathbb{Z}} {\mathbb{Z}}[{\pi}_1 X]
\cong K$,
in which the intersection form is hyperbolic.) Since $G$ is 
contained in $M$, all its double point loops are trivial in ${\pi}_1 N$. 
Cap $G$ by null-homotopies for its tips in $N$, to get
a collection of capped $(m+1)$-gropes $G^c$. These are not
capped gropes in the conventional sense, because the body has 
self-intersections (but these are ${\pi}_1 N$-null), and also
the caps may intersect any surface stage -- however the total number
of group elements represented by the double point loops of
$G^c$ is at most $m$, as observed above.

Split $G^c$ (as in Lemma \ref{splitting}), with respect to the group 
elements at its caps. In other words, first split the caps, separating 
different group elements, and then proceed down the grope, splitting 
surface stages into genus $1$ pieces. (When splitting the surface stages,
ignore their intersections with any other surfaces and caps.)
The result, for each $G^c_i$, is a capped $(m+1)$-grope with the base
surface of high genus, with all surfaces above the first stage
of genus $1$, and with each cap having double points (intersections
with other caps/surface stages) with just one group element.

Consider a genus $1$ piece of the base surface. It is a base of
a capped ``dyadic'' $(m+1)$-grope (all surface stages have genus $1$), 
so has $m+1$ caps. There are at most $m$ group
elements present at the caps, so two of the caps must have the
same group element. Contract the grope along these two caps,
and push off all other caps/surfaces intersecting them,
thus creating only ${\pi}_1$-null intersections.
This produces a collection of ${\pi}_1$-null  transverse
pairs of spheres, and reduces the problem to Chapter 6 of \cite{FQ}.
\qed


\begin{thebibliography}{10}
\setlength{\parskip}{.15cm}

\bibitem[D]{D} W. Dwyer, {\em Homology, Massey products and maps between
groups}, J. Pure Appl. Algebra 6 (1975), 177-190.

\bibitem[F]{F} M.H. Freedman, {\em The disk theorem for four-dimensional
manifolds}, Proc. ICM Warsaw (1983), 647-663.

\bibitem[FQ]{FQ} M.H. Freedman, F. Quinn, {\em The topology of 
4-manifolds}, Princeton Math. Series 39, Princeton, NJ, 1990.

\bibitem[FT]{FT} M.H. Freedman, P. Teichner, {\em $4$-Manifold Topology II:
Dwyer's filtration and surgery kernels}, 
Invent. Math. 122 (1995), 531-557.

\bibitem[K]{K} V.S. Krushkal, {\em Exponential separation in $4$-manifolds}, 
Geom. Topol. 4 (2000), 397-405.

\bibitem[KQ]{KQ} V.S. Krushkal, F.Quinn, {\em Subexponential groups in $4$-manifold 
topology}, Geom. Topol. 4 (2000), 407-430.

\end{thebibliography}
\end{document}